\documentclass[11pt]{amsart}

\usepackage{palatino}

\usepackage{latexsym}
\usepackage{amsmath,amsfonts,amssymb,epsfig}
\usepackage{amsthm}
\usepackage{hyperref}

\newtheorem{theorem}{Theorem}[section]
\newtheorem{definition}[theorem]{Definition}

\newtheorem{proposition}[theorem]{Proposition}
\newtheorem{lemma}[theorem]{Lemma}
\newtheorem{conjecture}[theorem]{Conjecture}

\newtheorem{rem}[theorem]{Remark}
\newenvironment{remark}{\begin{rem} \rm}{\end{rem}}
\newtheorem{ide}[theorem]{Idea}

\newtheorem{exa}[theorem]{Example}

\newtheorem{ques}[theorem]{Question}
\newenvironment{question}{\begin{ques} \rm}{\end{ques}}
\newtheorem{pro}[theorem]{Problem}

\newtheorem{spe}[theorem]{Speculation}

\newcommand{\Lie}{{\mathcal L}}

\font\bbb=msbm10 scaled 1100
\newcommand{\real}{\mbox{\bbb R}}       

\newcommand{\df}[1]{{\sc #1}}

\begin{document}

\title{OVERTWISTED ENERGY-MINIMIZING CURL EIGENFIELDS}

\author{Robert Ghrist}
\address{Department of Mathematics and Coordinated Sciences Laboratory,
        University of Illinois, Urbana IL, 61801}
\thanks{RG supported in part by NSF PECASE Grant \# DMS - 0337713.}

\author{Rafal Komendarczyk}
\address{Department of Mathematics, Georgia Institute of Technology,
        Atlanta GA, 30332-0160}
\thanks{RK supported in part by NSF Grant \# DMS - 0134408 .}


\begin{abstract}
This paper concerns topological and geometric properties of
energy-minimizing solutions to the steady Euler equations for a
fluid on a 3-dimensional manifold. Specifically, we consider
energy-minimizing divergence-free eigenfields of the curl operator
in dimension three from the perspective of contact topology. We give
a negative answer to a question of Etnyre and the first author by
constructing curl eigenfields which minimize $L^2$ energy on their
co-adjoint orbit, yet are orthogonal to an \df{overtwisted} contact
structure. We then give progress toward the conjecture that
$K$-contact structures on Seifert fibered manifolds always define
\df{tight} minimizers.
\end{abstract}

\maketitle

\section{Introduction}

Eigenfields of the curl operator form an important class of
solutions to the steady Euler equations in dimension three. These
equations model the velocity field of an inviscid, incompressible
fluid flow on a Riemannian manifold \cite{Khesin98}. It has been
observed \cite{Hamilton85,Chicone,Ghrist00_2} that there is a
correspondence between curl eigenfields and contact 1-forms in
dimension three. Let $M$ be a 3-manifold and $\Omega^1(M)$ denote
the 1-forms on $M$. Recall that a contact 1-form
$\alpha\in\Omega^1(M)$ is one which satisfies $\alpha\wedge
d\alpha\neq 0$. Such 1-forms have as their kernel a totally
nonintegrable plane field known as a {\em contact structure}. The
correspondence is this: any nonvanishing curl eigenfield is dual
to a contact 1-form via the Riemannian metric; and any contact
1-form can be realized as dual to a nonvanishing curl eigenfield
for some Riemannian structure.

This observation raises interesting questions concerning the
interplay between fluid dynamical properties of curl eigenfields and
topological properties of contact structures. More specifically, one
can investigate how the topological {tight/overtwisted} dichotomy
for contact structures relates to the physical properties of the
fluid like {energy}, {periodic orbits} etc. In \cite{Ghrist00_2} (p.
16) the authors asked if variational principles of the fluid flows
descend to variational principles for the corresponding contact
structures.

\begin{question}\label{qe:q1}
Does a nonvanishing curl eigenfield which minimizes ($L^2$) energy
on its coadjoint orbit under the volume-preserving diffeomorphism
group necessarily define a tight contact structure?
\end{question}

Using a recent result of the second author
\cite{Komendarczyk_nodal}, we give a negative answer to this
question. This involves constructing special $S^1$-invariant curl
eigenfields on products $S^1\times\Sigma$, where $\Sigma$ is a
closed orientable surface of genus $g>0$. In the second part of this
note, we demonstrate that under additional symmetry conditions, the
energy-minimization condition does yield a tightness constraint on
the associated contact structures.

Throughout the paper, we use the language of differential forms
and global analysis, as in \cite{Khesin98}. We restrict attention
to the class of volume-preserving vector fields and 1-forms
exclusively.


\section{Contact structures}

Let $(M,g)$ be a Riemannian 3-manifold. The \df{curl operator} on
1-forms is $\ast\,d:\Omega^1(M)\rightarrow\Omega^1(M)$, where $\ast$
is the Hodge star operator and $d$ the exterior derivative. An
\df{eigenform} of curl is any 1-form $\alpha$ satisfying
$\ast\,d\alpha=\mu\alpha$ for some $\mu\in\real$.

Given a nonvanishing curl eigenform $\alpha$ with nonzero eigenvalue
$\mu\neq 0$, it easily follows that
\[
    \alpha \wedge d\alpha = \alpha \wedge (\mu\ast\alpha) \neq 0.
\]
Therefore, $\alpha$ is a \df{contact form} and the plane field
$\xi=\text{ker}(\alpha)$ defines a {\it contact structure}
--- a nowhere integrable plane field on $M$. If $\alpha$ is the dual
1-form to a vector field, then $\xi$ is the orthogonal plane field.
All contact structures in this paper are transversely orientable, as
they are dual to globally-defined vector fields.

It has been known since the work of Bennequin and Eliashberg
\cite{Benn,Eli} that there are two fundamentally different classes
of contact structures.

\begin{definition}\label{def:tight_overt}
A contact structure $\xi$ is \df{overtwisted} if and only if there
exists an embedded disk $D^2\subset M$ such that $D$ is transverse
to $\xi$ near $\partial D$ but $\partial D$ is tangent to $\xi$. Any
contact structure which is not overtwisted is called \df{tight}.
\end{definition}

As one might expect from the definition, it is rather difficult to
determine if a given contact structure is tight or overtwisted. One
of the more successful recent techniques for solving the
classification problem involves examining the \df{characteristic
surface}.
\newpage
\begin{definition}
Let $X$ be a vector field preserving the contact plane distribution
$\xi$, i.e., $\Lie_X\xi=0$.
The characteristic surface $\Gamma_X$ is the set of tangencies of $X$
with $\xi$,
\begin{gather}\label{eq:char_surface}
    \Gamma_X=\{p\in M: X_p\in \xi_p\}.
\end{gather}
\end{definition}

The following result is essential for our study of $S^1$-invariant
curl eigenforms on circle bundles.

\begin{theorem}[Giroux \cite{Giroux01}]\label{th:giroux_bundle}
Let $\xi$ be an $S^1$-invariant contact structure on a principal
circle bundle $\pi:P\to\Sigma$ over a closed oriented surface
$\Sigma$. Let $\Gamma=\pi(\Gamma_{S^1})$ be a projection of the
characteristic surface $\Gamma_{S^1}$ onto $\Sigma$. Denote by
$e(P)$ be the Euler number of the bundle $P$. Then $\xi$ and all
covers of $\xi$ are tight if and only if one of the following holds:
\begin{itemize}
    \item[(i)] For $\Sigma\neq S^2$ none of the connected
        components of $\Sigma/\Gamma$ is a disc.
    \item[(ii)] For $\Sigma=S^2$, $e(P)<0$ and $\Gamma=\emptyset$.
    \item[(iii)] For $\Sigma=S^2$, $e(P)\ge 0$ and $\Gamma$ is connected.
   \end{itemize}
\end{theorem}
A contact structure all of whose covers are tight is called
\df{universally tight}.


\section{Energy and eigenvalues}

An important feature of any curl eigenform $\alpha$ is the fact that
it extremizes the $L^2$-energy, defined as,
\begin{gather}\label{eq:energy}
E(\alpha)=\|\alpha\|^2_{L^2}=\int_M \alpha\wedge\ast\alpha,
\end{gather}
among all 1-forms obtained from $\alpha$ by pullbacks through volume
preserving diffeomorphisms. The set of such forms, $\Xi_\alpha$, is
the \df{coadjoint orbit} of $\alpha$ under the action of the
volume-preserving diffeomorphism group of $M$:
\begin{gather}
    \Xi_\alpha = \{\beta: \beta=\varphi_\ast(\alpha),\ \varphi\in
    \text{Diff}_0(M),\ \varphi_\ast(\ast 1)=\ast 1\}.
\end{gather}

The question of energy minimization on the coadjoint orbit is more
delicate, and closely related to spectral data. The following result
is one of the few general results available:

\cite{Khesin98}:
\begin{proposition}\label{th:prop_minimizers}
A curl eigenform $\alpha_1$, (i.e. an eigenform of the curl operator
$\ast\,d:\mathcal{H}\to \mathcal{H}$,
$\mathcal{H}=\{\alpha:\delta\alpha=0\}=\{\text{``divergence free''
1-forms}\}$). corresponding to the first eigenvalue $\mu_1\neq 0$ is
a minimizer of the energy $E$ on $\Xi_{\alpha_1}$.
\end{proposition}
\begin{proof}
The operator $\ast\,d$ is elliptic and consequently it's analytic
realization is unbounded on $L^2_\mathcal{H}(M,\Lambda^1 T^\ast M)$
(an $L^2$- completion of $\mathcal{H}$, closed and self-adjoint; it
has a compact inverse $\ast\,d^{-1}$ defined on the orthogonal
complement of its kernel (see \cite{Yoshida90}).  We can also choose
an orthonormal basis of eigenforms $\{\alpha_i\}$ in
$(L^2_\mathcal{H}(M,\Lambda^1 T^\ast M),(\cdot,\cdot)_{L^2})$ such
that,
\begin{gather}\label{eq:eigenbeltrami}
\ast\,d^{-1}\alpha_i=\frac{1}{\mu_i}\,\alpha_i,\qquad 0<\mu^2_1\leq
\mu^2_2\leq\ldots\leq \mu^2_i\leq\ldots
\end{gather}
For an arbitrary $L^2$ 1-form $\beta\in \text{Im}(\delta)$ we have
\begin{gather}\label{eq:helicity_estim}
 |(\ast\,d^{-1}\beta,\beta)_{L^2}|=|\sum_i
\frac{1}{\mu_i}(\alpha_i,\beta)^2_{L^2}|
\leq\frac{1}{|\mu_1|}(\beta,\beta)_{L^2}=\frac{1}{|\mu_1|}E(\beta).
\end{gather}
One obtains a lower bound for the energy $E(\beta)$,
\begin{gather*}
E(\beta)\ge |\mu_1||(\ast\,d^{-1}\beta,\beta)_{L^2}|.
\end{gather*}
The above inequality becomes the equality if and only if $\beta$ is
a $\mu_1$- eigenform of $\ast\, d$. The claim follows from the fact
that the \df{helicity}, $(\ast\,d^{-1}\beta,\beta)_{L^2}$, is
invariant under volume preserving transformations see
\cite{Khesin98}.
\end{proof}

From now on we do not distinguish between operators defined on
various spaces of smooth differential forms and their analytic
realizations defined on $L^2$- completions of those spaces.

The curl operator $\ast\,d:\mathcal{H}\to \mathcal{H}$ is a
self-adjoint first-order elliptic operator, and the principal
eigenvalue $\mu_1$ enjoys a variational characterization through the
Rellich's quotient. Via Lemma \ref{th:ck_lemma} we have,
\begin{gather}
\label{eq:rellich_quot}    \mu_1 = \inf_{\alpha\in
\mathcal{H}^\perp_0}
    \frac{|(\ast\,d\alpha,\alpha)_{L^2}|}{\|\alpha\|^2_{L^2}}
    \quad
    \Leftrightarrow\quad \mu^2_1=\inf_{\alpha\in \mathcal{H}^\perp_0}
    \frac{(\Delta^1_M\alpha,\alpha)_{L^2}}{\|\alpha\|^2_{L^2}},
\\
\notag    \mathcal{H}^\perp_0=\text{Ker}(\ast\,d)^\perp=
    \{\alpha\in\Omega^1(M):\alpha=\delta\beta,
    \text{ for some }\beta\in\Omega^2(M)\},\
\end{gather}
Observe that the curl squared is equal to the Hodge Laplacian ,
$(\ast\,d)^2=\delta\,d$, on $\mathcal{H}$. Therefore any curl
eigenform $\alpha$, (i.e. $\ast\,d\alpha=\pm\mu\,\alpha$) is
automatically a co-closed $\mu^2$-eigenform of the Hodge Laplacian
i.e.
\begin{gather}\label{eq:betrami-lap}
\Delta^1_M\alpha=\delta\,d\,\alpha=\ast\,d\ast d\alpha=\mu^2\,\alpha,
\end{gather}
Clearly the curl $\ast\,d$ commutes with $\Delta^1_M$, therefore
both of these operators are simultaneously diagonizable on
$\mathcal{H}$ in a suitable orhonormal basis of curl eigenforms;
\begin{gather*}
\mathcal{H}=\bigoplus^\infty_{i=1} E^{\Delta} (\mu^2_i), \quad
E^{\Delta}(\mu^2_i)\perp E^{\Delta}(\mu^2_j),\ i\neq j
\qquad0<\mu^2_1\leq \mu^2_2\leq\ldots\leq \mu^2_i\leq\ldots
\end{gather*}
where $E^{\Delta}(\mu^2_i)$ stands for the $\mu^2_i$-eigenspace of
$\Delta^1_M$, and
\begin{gather*}
E^{\Delta}(\mu^2_i)=E^{\ast\,d}(\mu_i)\oplus E^{\ast\,d}(-\mu_i).
\end{gather*}
(We allow one of $E^{\ast\,d}(\mu_i)$, $E^{\ast\,d}(-\mu_i)$ to be
trivial.) We may conclude further that there exist two positive
operators $\sqrt{\Delta}_{+}$, $\sqrt{\Delta}_{-}$, such that
$\ast\,d=\sqrt{\Delta}_{+}-\sqrt{\Delta}_{-}$.

The following useful fact, which can be traced back to work in
\cite{Chandrasekhar57}, tells us how to effectively find a basis of
curl eigenforms from a basis of co-closed $\Delta^1_M$-eigenforms.
\begin{lemma}\label{th:ck_lemma}
Any curl $\mu$-eigenform is automatically a co-closed
$\mu^2$-eigenform of the Laplacian $\Delta^1_M$. Conversely, given a
co-closed $\mu^2$-eigenform $\alpha\in \Omega^1(M)$ of $\Delta^1_M$
there exists a corresponding $\pm\mu$-curl eigenform $\beta_\pm\in
\Omega^1(M)$ given by
\begin{gather}\label{eq:ck_beltrami}
\beta_\pm=\mu\,\alpha\pm\ast\,d\alpha
\end{gather}
\end{lemma}
\begin{proof}
The first claim follows from \eqref{eq:betrami-lap}. The second
claim we verify by a direct calculation. Let $\beta_\pm$ be defined
by \eqref{eq:ck_beltrami}, we will show: $\ast
d\,\beta_\pm=\pm\,\mu\,\beta_\pm$. Since $\ast d\ast d
=\delta\,d=\Delta^1_M|_{\mathcal{H}}$, and $\delta\alpha=0$ we
obtain
\begin{gather*}
\ast\,d\beta_\pm=\mu\ast d\alpha\pm\Delta^1_M\alpha.
\end{gather*}
Secondly  $\Delta^1_M\alpha=\mu^2\alpha$, therefore
\begin{gather*}
\ast\,d\beta_\pm=\mu\ast d\alpha\pm\mu^2 \alpha=\pm\mu\,\beta_\pm .
\end{gather*}
\end{proof}

\section{Overtwisted principal eigenfields}

For the purpose of producing an overtwisted principal curl
eigenform we assume that $M$ is a trivial bundle
$P=S^1\times\Sigma$ and the metric $g$ on $P$ is a product metric
$g=1\oplus g_\Sigma$ such that fibres are of constant length $l$.
We will construct our example using a sequence of lemmas.

In the case of the product $(P, g)$, $P=S^1\times\Sigma$, $g=1\oplus
g_\Sigma$, the space of smooth 1-forms $\Omega^1(P)$ decomposes with
respect to the $L^2$-inner product induced by the metric $g$ as,
\begin{gather}\label{eq:tan_norm}
\Omega^1(P)=\Omega^1_N(P)\oplus\Omega^1_T(P)
\end{gather}
where,
\begin{gather*}
\Omega^1_N(P)=\{\alpha\in \Omega^1(P): \alpha=f\eta,\ f\in C^\infty(P)\},
\\
\Omega^1_T(P)=\Omega^1_N(P)^\perp\cap\Omega^1(P),
\\
\Omega^1_T(P)=\{\alpha\in \Omega^1(P): \alpha(X_\eta)=0\},
\end{gather*}
with $\eta$ and $X_\eta$ being the tangent 1-form and vector field
(resp.) to the $S^1$- fibres of unit magnitude ($\|X_\eta\|_g=1$).
The following lemma is another consequence of the product metric
assumption.

\begin{lemma}\label{th:dec_lap}
The Laplacian $\Delta^1_P$ preserves $\Omega^1_N(P)$,
$\Omega^1_T(P)$ and for $\alpha=f\eta+\beta$, $f\eta\in
\Omega^1_N(P)$, $\beta \in \Omega^1_T(P)$ we have the following
formula for the Laplacian at a point $(t,q)\in S^1\times \Sigma$,
\begin{gather}\label{eq:lap-product}
    \Delta^1_P\alpha=(-\Lie^2_\eta f+\Delta^0_\Sigma f_t)\eta +
    (-\Lie^2_\eta\,\beta+\Delta^1_\Sigma\,\beta_t),\quad \text{at}\quad (t,q),
\end{gather}
where $f_t=f|_{\{t\}\times\Sigma}\in C^\infty(\Sigma)$ and
$\Delta^0_\Sigma$ is the scalar Laplacian on $\Sigma$. Similarly
$\beta_t=\beta|_{\{t\}\times\Sigma}$ and $\Delta^1_\Sigma$ is the
1-form Laplacian on $\Sigma$.
\end{lemma}

\begin{proof}
 The first claim follows immediately from the formula \eqref{eq:lap-product}. We
 justify \eqref{eq:lap-product} by a direct calculation in the
 $X_\eta$- invariant frame $\{X_1, X_2, X_3\}$, $X_1=X_\eta$, (denote the co-frame
 by $\{\eta_i\}$, $\eta_1=\eta$) on $P$, where $X_2$, $X_3$ are tangent to the $\Sigma$ fibers.
 Denote:
 $\nabla_i=\nabla_{X_i}$ and recall the following formulas (see e.g.
 \cite{Jost02})
\begin{eqnarray}
 \nabla_i X_j & = & \Gamma^k_{ij}\, X_k, \quad
 \nabla_i\eta_k= -\Gamma^k_{i j}\,\eta_j,\quad \Gamma^k_{ij}=-\Gamma^j_{ik}, \\
\label{eq:lem1_lap0}  \nabla^2_{i\,i}& =&
-\nabla_i\nabla_i+\Gamma^j_{i\,i}\nabla_j,\quad
\Delta^0=-\nabla^2_{i\,i}.
 \end{eqnarray}
 The well known Weitzenbl\"{o}ck formula (see \cite{Jost02} p. 138) for the
 $k$- from Laplacian $\Delta^k$ tells us $(k=1)$:
 \begin{eqnarray}
  \label{eq:lem1_lap1}\Delta^1 \alpha & = &-\nabla^2_{i\,i}\alpha-\eta_i\wedge (X_j\lrcorner\ R(X_i,
  X_j)\alpha),\quad \alpha\in \Omega^1(M),\\
  \notag R(X_i,X_j)\alpha & =
  &\nabla_i\nabla_j\alpha-\nabla_j\nabla_i\alpha-\nabla_{[X_i,
  X_j]}\alpha.
 \end{eqnarray}
In the product metric we may choose locally an $X_\eta$- invariant
frame $\{X_i\}$, meaning:
\begin{gather*}
[X_1,X_j]=-[X_j,X_1]=0,\quad [X_2,X_3]\in T\Sigma.
\end{gather*}
Consequently (by $\Gamma^k_{i\,j}=\frac{1}{2}\langle [X_i, X_j],
X_k\rangle-\langle [X_j, X_k], X_i\rangle+\langle [X_k, X_i],
X_j\rangle$):
 \begin{eqnarray}
\label{eq:lem1_gam1} \Gamma^k_{i\,j} & = & 0,\quad\text{if one of the indices }i,j,k=1 , \\
\label{eq:lem1_gam2} \nabla_1 \Gamma^k_{i\,j} & = & 0, \quad\text{for all }i,j,k\\
\label{eq:lem1_curv} R(X_i, X_j)\eta_r & = & 0, \quad\text{if one of
the indices $i,j,k,r=1$, }
 \end{eqnarray}
 where \eqref{eq:lem1_curv} is a consequence of the following
 \begin{gather*}
  X_k\lrcorner\, R(X_i, X_j)\eta_r=\nabla_j\Gamma^r_{i\,k}-\nabla_i\Gamma^r_{j\,
  k}+\Gamma^r_{j\,n}\,\Gamma^n_{i\,k}-\Gamma^r_{i\,n}\Gamma^n_{j\,k}+(\Gamma^n_{i\,j}-\Gamma^n_{j\,i})\Gamma^r_{n\,k}.
 \end{gather*}
 In turn we obtain,
 \begin{gather}\label{eq:lem1_covfrm}
\nabla\eta_1=0,\quad \nabla_1\eta_k=0,\quad
\nabla_2\eta_i=-\Gamma^2_{i\,3}\,\eta_3,\quad
\nabla_3\eta_i=-\Gamma^3_{i\,2}\,\eta_2\ .
 \end{gather}
 Let $\alpha=f\eta+\beta$, then
 $\Delta^1_P\alpha=\Delta^1_P(f\eta_1)+\Delta^1_P\beta$, we obtain
 from \eqref{eq:lem1_lap1}, \eqref{eq:lem1_curv},
\begin{eqnarray*}
\Delta^1_P(f\eta_1) & = &
-\nabla^2_{i\,i}(f\eta_1)=-\nabla_i\nabla_i(f\eta_1)+\Gamma^j_{i\,i}\nabla_j(f\eta_1)\\
\eqref{eq:lem1_covfrm}: & = & (-\nabla_i\nabla_i f+\Gamma^j_{i\,i}\nabla_j f)\eta_1\\
\text{\eqref{eq:lem1_gam1}\ -\ \eqref{eq:lem1_curv}}: & = &
(-\nabla_1\nabla_1 f-\nabla_2\nabla_2 f+\Gamma^3_{2\,2}\nabla_3
f-\nabla_3\nabla_3 f+\Gamma^2_{3\,3}\nabla_2 f)\eta_1
 \end{eqnarray*}
 Treating $f=f_t$ as a family of functions $f_t\in C^\infty(\Sigma)$ dependent on $t\in S^1$, by \eqref{eq:lem1_lap0} we
 have,
 \begin{gather*}
-\nabla_2\nabla_2 f_t+\Gamma^3_{2\,2}\nabla_3 f_t-\nabla_3\nabla_3
f_t+\Gamma^2_{3\,3}\nabla_2 f_t=\Delta^0_\Sigma f_t,\quad
\text{therefore}\\
\Delta^1_P(f\eta)=(-\nabla_1\nabla_1 f + \Delta^0_\Sigma
f_t)\eta=(-\Lie^2_\eta f + \Delta^0_\Sigma f_t)\eta.
 \end{gather*}
Similar reasoning applies to $\beta=a_2\,\eta_2+a_3\,\eta_3$, by
\eqref{eq:lem1_covfrm} and \eqref{eq:lem1_curv} we have
\begin{gather*}
(-\nabla_1\nabla_1+\Gamma^j_{1\, 1}\nabla_j)\beta-\eta_1\wedge
(X_j\lrcorner\ R(X_1,
  X_j)\beta)=
-\nabla_1\nabla_1(a_2\eta_2+a_3\eta_3) \\  =-(\nabla_1\nabla_1
a_2)\eta_2-(\nabla_1\nabla_1 a_3)\eta_3=-\Lie^2_\eta\beta,
\end{gather*}
Treating $\beta=\beta_t$ as a family of 1-forms $\beta_t\in
\Omega^1(\Sigma)$, and using \eqref{eq:lem1_covfrm},
\eqref{eq:lem1_curv} and \eqref{eq:lem1_lap1} one shows,
\begin{gather*}
\sum^3_{i=2}\{-\nabla_i\nabla_i\beta+\Gamma^j_{i\,
i}\nabla_j)\beta-\eta_i\wedge (X_j\lrcorner\ R(X_i,
  X_j\beta)\}=\Delta^1_\Sigma \beta_t.
\end{gather*}
\end{proof}

\begin{lemma}\label{th:eigenvalues_lemma}
On the product manifold $P=S^1\times \Sigma$, $g=1\oplus g_\Sigma$
with constant length $l$ fibres. The first eigenvalue $\mu_1$ of the
curl operator satisfies,
\begin{gather*}
    \mu^2_1=\min\left\{\nu_1,\Bigl(\frac{2\pi}{l}\Bigr)^2\right\},
    \quad\text{where}\quad
    \nu_1=\inf_{f\in L^2(\Sigma),\
    f\neq\text{const}}\left\{\frac{(\Delta^0_\Sigma
    f,f)_{L^2}}{\|f\|^2_{L^2}}\right\}.
\end{gather*}
\end{lemma}

\begin{proof}
From the decomposition \eqref{eq:tan_norm} and the fact that
$\Delta^1_P$ preserves $\Omega^1_T(P)$ and $\Omega^1_N(P)$ (see
Lemma \ref{th:dec_lap}) we have
\begin{gather}\label{eq:eigenvalues}
    \mu^2_1=\min\left\{\mu^2_{1,T}, \mu^2_{1,N}\right\};\quad
    \mu^2_{1,r}=\inf_{\alpha\in \mathcal{H}^\perp_0\cap
    \Omega^1_r(P)}\left\{\frac{(\Delta^1_P\alpha,\alpha)_{L^2}}{\|\alpha\|^2_{L^2}}\right\}\
    \ r=T,N.
\end{gather}
In order to calculate $\mu^2_{1,N}$ notice that for any $\alpha\in
\mathcal{H}\cap\Omega^1_N(P)$, $\alpha=f\,\eta$, the function $f$
is constant on the fibres; hence $f\in C^\infty(\Sigma)$. Indeed,
$\delta\alpha=0$, and, since $\nabla\eta=0$ in the adapted frame
$\{X_1,X_2,X_3\}$ with $X_1=X_\eta$, we obtain
\begin{gather*}
0=\delta\alpha=X_i\lrcorner\nabla_i\alpha=X_i\lrcorner(\nabla_i
f\eta+f \nabla_i\eta)=\nabla_1 f=X_\eta f.
\end{gather*}
From \eqref{eq:lap-product} we conclude that for any
$\alpha=f\,\eta\in\mathcal{H}\cap\Omega^1_N(P)$, we have:
$\Delta^1_\Sigma \alpha=(\Delta^0_\Sigma f)\eta$ and consequently,
\begin{gather*}
    (\Delta^1_P\alpha,\alpha)_{L^2}=(\Delta^0_\Sigma f\eta,f\eta)_{L^2}
    =\int_{S^1\times\Sigma} (f\,\Delta^0_\Sigma
    f)\eta\wedge\ast\eta\\
    =\int_{S^1}\int_\Sigma f\,\Delta^0_\Sigma f =l\,\int_\Sigma f\,\Delta^0_\Sigma f=l\,(\Delta^0_\Sigma
    f,f)_{L^2},\\
    \quad\text{ and }\quad\|\alpha\|^2_{L^2}=(\alpha,\alpha)_{L^2}=\int_{S^1\times\Sigma}
    f^2\eta\wedge\ast\eta=l\,\|f\|^2_{L^2}
\end{gather*}
where $\eta\wedge\ast\eta=\ast 1$. As a result,
\eqref{eq:rellich_quot},
\begin{gather}\label{eq:eigenvalue_N}
    \mu^2_{1,N}=\nu_1,\qquad \nu_1=\inf_{f\in L^2(\Sigma),\
    f\neq\text{const}}\left\{\frac{(\Delta^0_\Sigma
    f,f)_{L^2}}{\|f\|^2_{L^2}}\right\}.
\end{gather}
In other words $\mu^2_{1,N}$ is equal to the first eigenvalue of the
scalar Laplacian $\Delta^0_\Sigma$ on $\Sigma$.

In order to calculate $\mu^2_{1,T}$ we first calculate the
orthogonal basis of eigenforms of $\mathcal{H}\cap\Omega^1_T(P)$.
Let $\{\beta_m\}$ be an orthonormal basis of $\Delta^1_\Sigma$-
eigenforms on $\text{Ker}(\Delta^1_\Sigma)\oplus
\text{Im}(\delta)\subset L^2(\Lambda^1 T^\ast \Sigma)$, define for
all $m, n\in \mathbb{Z}^+$:
\begin{gather}\label{eq:eigen-basis}
h_0=g_0=1,\quad h_n=\cos\bigl(\frac{2\pi n t}{l}\bigr),\ g_n=\sin\bigl(\frac{2\pi n
t}{l}\bigr),
\\
\notag \alpha^g_{n\,m}=g_n\,\beta_m,\qquad\text{ and }
\qquad\alpha^h_{n\,m}=h_n\,\beta_m .
\end{gather}
Clearly $\{\alpha^g_{n\,m},\alpha^h_{n\,m}\}$ is a set of
1-eigenforms of $\Delta^1_P$ on $\mathcal{H}\cap\Omega^1_T(P)$.
Indeed, from \eqref{eq:lap-product} we have,
\begin{gather*}
\Delta^1_P\alpha^r_{mn}=\gamma^r_{mn}\alpha^r_{mn},\quad
\gamma^r_{mn}=\Bigl(\frac{2\pi n}{l}\Bigr)^2+\tilde{\nu}_m
\end{gather*}
where $\tilde{\nu}_m$ is the $m$-th eigenvalue of $\Delta^1_\Sigma$.
One easily shows that $\{\alpha^g_{n\,m},\alpha^h_{n\,m}\}$ is an
orthonormal basis of $\mathcal{H}\cap\Omega^1_T(P)$.
Consequently, all eigenforms of $\Delta^1_P$ on
$\mathcal{H}\cap\Omega^1_T(P)$ are listed in \eqref{eq:eigen-basis},
and we have
\begin{gather}\label{eq:eigenvalue_T}
\mu^2_{1,T}=\min\left\{\Bigl(\frac{2\pi}{l}\Bigr)^2,\tilde{\nu}_1\right\}
\end{gather}
It is left to show that $\tilde{\nu}_1=\nu_1$. By the Hodge
decomposition theorem:
\begin{gather*}
\Omega^1(\Sigma)=
\text{Ker}(\Delta^1_\Sigma)\oplus\text{Im}(d_\Sigma)
\oplus\text{Im}(\delta_\Sigma).
\end{gather*}
Moreover, $\Omega^0(\Sigma)\simeq\Omega^2(\Sigma)$ through the
Hodge-star isometry. We conclude that
$\text{Im}(\delta_\Sigma)=\{\ast_\Sigma\,d\, f;f\in
C^\infty(\Sigma)\}$. Since $\Delta^1_\Sigma$ commutes with
$\ast_\Sigma\,d$, any $\nu_m$-eigenfunction $f_m$ results in a
$\nu_m$-eigenform $\ast_\Sigma\,d\, f_m$. Therefore
$\tilde{\nu}_1=\nu_1$ and the lemma follows from
\eqref{eq:eigenvalues}, \eqref{eq:eigenvalue_N}, and
\eqref{eq:eigenvalue_T}.
\end{proof}

In \cite{Komendarczyk_nodal} the following theorem was proved:

\begin{theorem}[\cite{Komendarczyk_nodal}]\label{th:payne}
For an arbitrary closed compact orientable surface $\Sigma$, there
exists a smooth metric $g_\Sigma$ such that a nodal set
$f^{-1}_1(0)$ of the principal eigenfunction $f_1$ of
$\Delta^0_\Sigma$ is a single embedded circle which bounds a disc in
$\Sigma$.
\end{theorem}

Combining the above theorem with Lemma \ref{th:eigenvalues_lemma}
results in the following,

\begin{theorem}\label{th:principal_overt}
Let $\Sigma\neq S^2$ be an orientable surface of an arbitrary
nonzero genus. One can prescribe a metric $g_\Sigma$ on $\Sigma$
such that there exist an overtwisted curl eigenfield $v$ on the
product manifold $(S^1\times\Sigma,1\oplus g_\Sigma)$ which
minimizes the energy \eqref{eq:energy} on the coadjoint orbit
$\Xi_\alpha$.
\end{theorem}

\begin{proof}
In the first step we choose a metric $g_\Sigma$ on $\Sigma\neq S^2$
constructed in Theorem \ref{th:payne} and assume that the length of
fibres in $(S^1\times\Sigma,1\oplus g_\Sigma)$ is equal to $l$. By
Lemma \ref{th:eigenvalues_lemma} we may choose $l$ small so that the
first eigenvalue satisfies $\mu_1=\nu_1$. The proof of
\ref{th:eigenvalues_lemma} implies that the corresponding eigenspace
$E^\Delta(\mu^2_1)$ is spanned by two independent co-closed
$\mu^2_1$- eigenforms of $\Delta^1_P$: $\alpha_1=f_1\,\eta$, and
$\alpha_2=\ast\,d_\Sigma f_1$. (The dimension is two since
$g_\Sigma$ can be chosen in a residual subset of an open set of
metrics: see, \cite{Komendarczyk_nodal}.) By earlier considerations
$E^\Delta(\mu^2_1)=E^{\ast\,d}(\mu_1)\oplus E^{\ast\,d}(-\mu_1)$ and
$E^\Delta(\mu^2_1)$ is spanned by two independent $\pm\mu_1$- curl
eigenforms. Choosing any linear combination of $\alpha_1$ and
$\alpha_2$ Lemma \ref{th:ck_lemma} leads to $\pm\mu_1$- curl
eigenforms given by,
\begin{gather*}
\beta_{\pm}=f_1\eta\pm\ast_{\Sigma}\,d\,f_1.
\end{gather*}

These forms are nonvanishing since the set of zeros is clearly equal
to the singular part of the nodal set of $f_1$. This singular part
is empty for a generic choice of metric (see \cite{Uhlenbeck76}).
Both forms are $S^1$-invariant and overtwisted by Theorem
\ref{th:giroux_bundle}. Indeed, the projection of the characteristic
surface $\Gamma_{S^1}$ of $\alpha_{\pm}$ onto $\Sigma$ is equal to
$f^{-1}_1(0)$, the nodal set of $f_1$. By the choice of the metric
$\pi(\Gamma_{S^1})$ bounds a disk. Now, the dual curl eigenfields
$\beta^{\#}_{\pm}$ minimize energy \eqref{eq:energy} on
$\Xi_{\beta_{\pm}}$ due to Proposition \ref{th:prop_minimizers}.
\end{proof}
\begin{remark}\label{rm:perturb}
 By perturbing a product metric on $P$ the eigenvalues $\pm\mu_1$,
 ``split apart'' giving only one of minimal absolute eigenvalue, and one dimensional eigenspaces. If
 the perturbation is small then the resulting eigenform will be
 $C^0$- close to $\beta_+$( or $\beta_-$), and define an isotopic
 contact structure (by Grey's Theorem, \cite{Eli}), which in turn
 must be an overtwisted minimizer.
\end{remark}


\section{Symmetry and tight energy minima}

There are certain cases for which principal curl eigenforms can
never be overtwisted. Courant's theorem on nodal sets quickly yields
such a rigidity result for $S^1\times S^2$:

\begin{proposition}
Let $M=S^1\times S^2$ with any product metric giving the $S^1$
fibres a constant length $l$. Then the principal eigenform of curl
on $M$ is never overtwisted.
\end{proposition}
\begin{proof}
It follows from Lemma \ref{th:eigenvalues_lemma} that if the first
eigenvalue $\mu_1$ is equal to $\frac{2\pi}{l}$ then the curl
eigenforms have zeros (since the dual vector fields are tangent to
$S^2$). If $\mu_1=\nu_1$ Theorem \ref{th:giroux_bundle} tells us
that a contact form is tight if and only if the projection of
$\Gamma_{S^1}=f_1^{-1}(0)$ is a single circle. On the other hand
Courant's theorem implies that the principal eigenfunction on a
closed surface always has exactly two nodal domains, which in turn
implies tightness.
\end{proof}

For the more general case, tightness can be forced by additional
symmetry. Let $(M,g_M)$ be a Riemannian 3- manifold which admits a
unit Killing vector field $X$ orthogonal to a contact structure
$\xi=X^\perp$. Choosing a local frame of vector fields
$\{X_1,X_2,X_3\}$, $X_1=X$, since $\Lie_X g_P=0$ and consequently:
$g_P(\nabla_V X, W)=-g_P(V, \nabla_W X)$, we obtain
\begin{gather}
\label{eq:kcont_christ1}\Gamma^2_{1\,1}=\Gamma^3_{1\,1}=\Gamma^2_{2\,1}=\Gamma^3_{3\,1}=0,\qquad
\Gamma^k_{i\,j}=-\Gamma^j_{i\,k},\\
\label{eq:kcont_christ2}\Gamma^2_{3\,1}=\Gamma^3_{2\,1}=\lambda,
\end{gather}
where $\nabla_i X_j=\Gamma^k_{i\,j} X_k$, and $\lambda$ is an
invariant independent on the choice of $\{X_2, X_3\}$. Cartan's
equations of structure imply that dual 1-form $\eta=X^\flat$
satisfies
\begin{gather}
\ast\,d\,\eta=2 \lambda\,\eta,\qquad \lambda\in C^\infty(P).
\end{gather}
In the case $\lambda(x)=\lambda=\text{const}$, the triple
$(M,\eta,g_M)$ is called \df{$(K,\lambda)$- structure} on $M$
(\cite{Nicolaescu98}), in which case $\eta$ defines a curl eigenform
on $M$. The $(K,1)$- manifolds are also known as \df{$K$-contact}
manifolds or equivalently, in dimension 3, Sasakian manifolds,
(\cite{Blair02, Nicolaescu98}), the Hopf fields on the round $S^3$
are the classical example (see \cite{Khesin98}). We remark here that
if $M$ admits a $(K,\lambda)$- structure then it is topologically a
Seifert manifold (see \cite{Nicolaescu98} for the proof).

Now we focus on the special case of $(M,\eta,g_M)$, namely the case
of a principal $S^1$-bundle over a closed orientable surface
$\Sigma$, equipped with a bundle metric $g_P$. All fibres are
geodesics, as the bundle metric is invariant under the action of a
unit Killing field $X$ tangent to the fibres, and the projection
$\pi:P\rightarrow \Sigma$ is a Riemannian submersion.

By Theorem \ref{th:giroux_bundle} $\eta$ is necessarily tight since
the contact plane distribution is orthogonal to the fibres.

Denote by $\mathcal{H}_{S^1}$ the subspace of $S^1$- invariant
1-forms in $\mathcal{H}\subset\Omega^1(P)$.
\begin{proposition}\label{th:K-cont_min}
Any curl eigenform $\eta$ defined by a $(K,\lambda)$-structure on
$(P,g_P)$ is always energy-minimizing on
$\mathcal{H}_{S^1}\cap\Xi_\eta$. Let $\nu$ be the first nonzero
eigenvalue of the scalar Laplacian $\Delta^0_\Sigma$ on $\Sigma$. If
$\nu>3\lambda^2$ then $\eta$ is a principal curl eigenform on
$\mathcal{H}_{S^1}$.
\end{proposition}

\begin{proof}
We provide the proof of the first claim for $\lambda>0$ (in the case
of $\lambda<0$ the reasoning is analogous). The space
$\mathcal{H}_{S^1}$ decomposes as
\begin{gather*}
\mathcal{H}_{S^1}=\mathcal{H}^+_{S^1}\oplus\mathcal{H}^-_{S^1},
\end{gather*}
where $\mathcal{H}^\pm_{S^1}$ is a subspace spanned by
positive/negative curl eigenforms. We need to show that $\eta$ is an
energy minimizer on $\mathcal{H}_{S^1}\cap\Xi_\eta$. Given a volume
preserving diffeomorphism $\varphi:P\rightarrow P$ we denote
$\eta_\varphi=\varphi_\ast(\eta)\in \Xi_\eta$. Under the assumptions
on the $\varphi$ action, $\eta_\varphi\in\mathcal{H}_{S^1}\cap\Xi_\eta$.
We expand $\eta_\varphi$ in the eigenbasis of
curl eigenforms \eqref{eq:eigenbeltrami}, $\eta_\varphi=\sum_{i\ge
0} c^{+}_i\alpha^{+}_i+\sum_{i< 0} c^{-}_i\alpha^{-}_i$, where
$\{\alpha^{\pm}_i\}$ span $\mathcal{H}_{S^1}$. Since the
helicity $(\ast d^{-1}\eta_\varphi,\eta_\varphi)$ is invariant under
$\varphi$, as in \eqref{eq:helicity_estim}, we obtain
\begin{gather*}
0<\frac{E(\eta)}{2\lambda}=(\ast d^{-1}\eta,\eta)=(\ast
d^{-1}\eta_\varphi,\eta_\varphi)=\sum_{i\ge 0}
\frac{(c^{+})^2}{\mu^{+}_i}+\sum_{i< 0} \frac{(c^{-})^2}{\mu^{-}_i}
\end{gather*}
where $\mu^\pm_i$, positive/negative eigenvalues of $\ast d$ on
$\mathcal{H}_{S^1}$. Since the second sum is negative we can
estimate $\mu^+_1 (\ast d^{-1}\eta,\eta)\leq E(\eta_\varphi)$. To
finish the proof it suffices to show that $2\lambda=\mu^+_1$, then we
obtain $E(\eta)\leq E(\eta_\varphi)$ which proves the claim.

We derive the equality $2\lambda=\mu^+_1$ by a calculation in a
local (co)frame $\{X_i\}$($\{\eta_i\}$), $X_1=X,$, where $X$ is the
Killing field tangent to the fibres and
$\ker{\alpha_1}=\text{span}\{X_2,X_3\}$. We adapt the notation from
Lemma \ref{th:ck_lemma}. Let $\alpha_1=a_i\eta^i=f\eta+\beta$ be the
curl eigenform satisfying
\begin{gather}\label{eq:alpha_1}
 \ast d\,\alpha_1=\mu^+_1\alpha_1.
\end{gather}
Since, $\Lie_X\alpha_1=0$, using the Cartan formula we obtain:
 \begin{eqnarray*}
 0 & = & \Lie_{X} \alpha = \mu^+_1\,X\lrcorner\ast\alpha + df\\
-df=-\nabla_i f\,\eta^i & = & \mu^+_1\,X_1\lrcorner\ast\alpha\\
df=\nabla_1 f\,\eta^1+\nabla_2 f\,\eta^2+\nabla_3 f\,\eta^3 & = &
\mu^+_1\,(a_2\,\eta^3-a_3\,\eta^2).
 \end{eqnarray*}
This leads to the following equations,
 \begin{gather}\label{eq:eq_Lie}
    \nabla_1 f=0,\quad
    \nabla_2 f=-\mu^+_1\,a_3,\quad
    \nabla_3 f=\mu^+_1\,a_2.
\end{gather}
Applying \eqref{eq:alpha_1} and \eqref{eq:eq_Lie} to
$d\,\alpha=\nabla_i a_k\eta^i\wedge\eta^k -
a_k\Gamma^k_{ji}\,\eta^i\wedge\eta^j$ yields
\begin{gather}\label{eq:der1}
\mu^+_1\, f=\frac{1}{\mu^+_1}(-\nabla_2\nabla_2 f-\nabla_3\nabla_3
f+\Gamma^2_{33}\nabla_2 f+\Gamma^3_{22}\nabla_3 f)+ 2\lambda f.
\end{gather}
By \eqref{eq:lem1_lap0} and \eqref{eq:kcont_christ1},
\eqref{eq:kcont_christ2} we obtain the following equation for $f$:
\begin{gather}\label{eq:scalar_lap}
 \Delta^0_{P}f=\mu^+_1(\mu^+_1-2\lambda)f.
\end{gather}
Equations \eqref{eq:eq_Lie} imply that for $\alpha_1$ to be
nontrivial, $f$ cannot be a constant zero function. Since
$\Delta^0_{P}$ is a positive operator we conclude that
$\mu^+_1\geq 2\lambda$; consequently, $\mu^+_1=2\lambda$ because
$\mu^+_1$ is the first positive eigenvalue.

Equations \eqref{eq:der1} and \eqref{eq:scalar_lap} are valid for
any $S^1$-invariant $\mu$-eigenform. Because $\pi$ is the Riemannian
submersion: $\Delta^0_P(h\circ\pi)=\pi\circ\Delta^0_\Sigma h$, $h\in
C^\infty(\Sigma)$, and the proof of the second statement follows
from the equation: $\nu=\mu(\mu-2\lambda)$. Indeed, for
$\gamma^2=\nu$, $\gamma>0$, we obtain
\begin{gather}
 \mu^2-2\,\lambda\,\mu+\gamma^2=0,\qquad
 \overline{\mu}=\lambda+\sqrt{\lambda^2+\gamma^2},\
 \underline{\mu}=\lambda-\sqrt{\lambda^2+\gamma^2},
\end{gather}
where $\overline{\mu}$, $\underline{\mu}$ are the roots of the
equation. Consequently,
\begin{gather*}
\mu_1^-\in(-\infty,-\gamma),\ \text{if}\quad \mu_1^+=\lambda>0\\
\mu_1^+\in(\gamma,+\infty),\ \text{if}\quad \mu_1^-=\lambda<0,
\end{gather*}
and it suffices to assume $\nu=\gamma^2>3\lambda^2$, in order to
assure $\lambda$ to be the principal eigenvalue of $\ast\,d$ on
$\mathcal{H}_{S^1}$.
\end{proof}

We may think about $\lambda\neq 0$ as a ``topological deviation''
from the $\lambda=0$ case. We note that Hopf fields are principal
curl eigenfields of $\ast\,d$ on $\mathcal{H}$ and therefore, by
Proposition \ref{th:prop_minimizers}, energy minimizers. Lemma
\ref{th:eigenvalues_lemma} and \ref{th:ck_lemma}, provide a
methodology to construct energy minimizers on products: $S^1\times
\Sigma$. In the case of a $(K,\lambda)$- structure, however,  $\eta$
seem to be the most natural candidate for the energy minimizer on
$\Xi_\eta$. The following special case, shows when $\eta$ becomes
the principal curl eigenform on $\mathcal{H}$.

We choose a special bundle metric $\hat{g}_P$ on a principle circle
bundle $P$, of constant length $l$, $S^1$- fibers. Let $X$ a unit
Killing vector field, tangent to the fibers such that $\eta=X^\flat$
defines a $(K,\lambda)$- structure on $P$. In addition we assume
that there always exists an $X$- invariant local frame of vector
fields $\{X_1,X_2,X_3\}$, $X_1=X$, on $P$, namely
\begin{gather}\label{eq:zero_brackets}
[X_1,X_i]=0,\qquad i=1,2,3.
\end{gather}
Under this assumption and \eqref{eq:kcont_christ1},
\eqref{eq:kcont_christ2} we obtain
\begin{gather}\label{eq:pres_lie_frms}
\Lie_X\eta_i=0,\qquad \eta_i=X^\flat_i.
\end{gather}

\begin{theorem}\label{th:princ_k-lamb}
The curl eigenform $\eta$ defined by a $(K,\lambda)$-structure on
$(P,\hat{g}_P)$ is a principal curl eigenform on $\mathcal{H}$, if
\begin{gather}
\lambda^2<\min\bigl(\frac{\nu}{3},\frac{4\pi^2}{l^2}\bigr),
\end{gather}
where $\nu$ is the first nonzero eigenvalue of the scalar Laplacian
$\Delta^0_\Sigma$ on $\Sigma$, and $l$ the length of the fibre.
\end{theorem}
\begin{proof}
 Observe that on $(P,\hat{g}_P)$ the operator $-\Lie^2_X$ commutes with
 $\Delta^1_P=\delta\,d+d\,\delta$. We check that $\Lie_X:\Omega^\ast(P)\to \Omega^\ast(P)$
 commutes with the Hodge star operator $\ast: \Omega^\ast(P)\to
 \Omega^\ast(P)$. Indeed since $\Lie_X$ respects the wedge product:
 \begin{gather*}
  \Lie_X (\omega_1\wedge\omega_2)=\Lie_X \omega_1\wedge\omega_2+\omega_1\wedge\Lie_X
  \omega_2,\quad \omega_1\in \Omega^j(P),\omega_2\in \Omega^k(P),
 \end{gather*}
 and the property \eqref{eq:pres_lie_frms}, for any $k$- form
 $\alpha=\sum_I a_I\omega_I$, where  $\omega_I=\eta_{i_1}\wedge \eta_{i_2}
 \wedge\ldots\wedge \eta_{i_k}$, $I=(i_1,\ldots,i_k)$ we have
 \begin{gather*}
 \Lie_X \alpha= \sum_I (\Lie_X a_I)\omega_I,\qquad \text{since }\Lie_X\omega_I=0.
 \end{gather*}
Consequently, $\ast\,\Lie_X=\Lie_X\,\ast$ follows from the
definition of the Hodge star operator,
\begin{gather}
\ast\Lie_X\alpha=\sum_I (\Lie_X a_I)\ast\omega_I=\Lie_X \sum_I a_I
\ast\omega_I=\Lie_X \ast\alpha,
\end{gather}
since $\Lie_X\omega_I=\Lie_X\ast\omega_I=0$. Hence, we obtain
$\Lie_X\Delta^1_P=\Delta^1_P\,\Lie_X$, because $\Lie_X$ commutes
with an exterior derivative $d$ and $\delta=\pm\ast d\ast$. This in
turn implies
\begin{gather}\label{eq:comlie_lap}
 (-\Lie^2_X)\Delta^1_P=\Delta^1_P(-\Lie^2_X).
\end{gather}
Consequently, we may define the decomposition
\begin{gather}\label{eq:decomp_lap}
\Delta^1_P=-\Lie^2_X+\Delta^H_P.
\end{gather}
where $\Delta^H_P=\Delta^1_P+\Lie^2_X$, and we call $\Delta^H_P$ a
horizontal Laplacian. Both $-\Lie^2_X$ and $\Delta^H_P$ are not
elliptic since take into account only derivatives in certain
directions, but $-\Lie^2_X$, and $\Delta^H_P$, have discrete
spectra, and  commute by \eqref{eq:comlie_lap}. Consequently, these
operators are simultaneously diagonalizable in a suitable $L^2$-
orhonormal basis of $\Delta^1_P$- eigenforms (see also
\cite{Bergery}). The decomposition \eqref{eq:decomp_lap} and
\eqref{eq:comlie_lap} imply that any eigenvalue $\tau$ of
$\Delta^1_P$ is a sum $\tau=\psi+\mu$ of eigenvalue: $\psi$ of
$-\Lie^2_X$, and $\mu$ of $\Delta^H_P$.

If $-\Lie^2_X \alpha=\psi\,\alpha$, $\alpha=a_i\,\eta_i$, then
\begin{gather*}
-\Lie^2_X a_i=\psi\, a_i.
\end{gather*}
Solving this equation in a local chart $(t,q)\in U\simeq S^1\times
V\subset P$, $V\subset\Sigma$, gives us $-\Lie^2_X a_i=-\partial^2_t
a_i(t,q)=\psi\, a_i(t,q)$, for a fixed $q$:
$a_i(t,q)=A_q\,\cos(\psi\, t)+B_q\,\sin(\psi\, t)$. Therefore
$\psi=\bigl(2\,\pi\,n/l\bigr)^2$, $n\in\mathbb{N}$. Now the theorem
is a consequence of: $\mathcal{H}_{S^1}=\text{Ker}(-\Lie^2_X)\cap
\mathcal{H}$, and Proposition \ref{th:K-cont_min}.
\end{proof}

The $(K,\lambda)$-structures of Theorem \ref{th:princ_k-lamb}, (i.e.
satisfying the assumption \eqref{eq:zero_brackets}) are completely
described in \cite{Nicolaescu98}. Their construction is based on the
Boothby-Wang construction (see also \cite{Blair02}). One builds the
metric $\hat{g}_P$ on $P$ as follows: choose a Riemannian metric
$g_\Sigma$ on $\Sigma$ and a connection $H_\omega$ on $P$ such that
the connection 1-form $\omega$ satisfies,
\begin{gather}\label{eq:conn}
d\omega=-\Omega=-2\,\mathcal{E}\, \sigma,\quad \sigma=\pi_\ast
(*_\Sigma 1),\quad \pi:P\to\Sigma,
\end{gather}
where $\mathcal{E}$ is a degree of the bundle $P$, and $\Omega$ the
curvature form of the connection. We define $\eta=\omega/\gamma$,
$\gamma=\text{const}$, and $X$ to be the unique vertical vector
field such that $\eta(X)=1$. Now we define $\hat{g}_P$ as follows,
\begin{gather*}
\hat{g}_P(X,X)=1,\quad \hat{g}_P(V,W)=g_\Sigma(\pi^\ast V,\pi^\ast
W),\quad \text{if $V, W$ horizontal}.
\end{gather*}
Clearly $\Lie_X \hat{g}_P=0$, therefore $X$ is a Killing vector
field, and by \eqref{eq:conn},
\begin{gather*}
d\omega=\gamma\,d\,\eta=-2\,\mathcal{E}
\sigma=-2\,\mathcal{E}\ast\eta,
\end{gather*}
which results in $\ast\,d\eta=(-2\mathcal{E}/\gamma)\eta$, and a
$(K, -\mathcal{E}/\gamma)$ structure: $(P,\eta,\hat{g}_P)$. By this
construction $\hat{g}_P$ also satisfies \eqref{eq:zero_brackets}.
\begin{conjecture}
The curl eigenform $\eta$ defined by a $(K,\lambda)$-structure is
always a tight energy-minimizer on $\Xi_\eta$.
\end{conjecture}

\small

\section{Acknowledgments}
We would like to thank the second referee for insightful comments,
especially for Remark \ref{rm:perturb}. The second author is
grateful to Jason Cantarella for the inspiring discussion during
Georgia Topology Conference in Athens 2000.


\begin{thebibliography}{10}

\bibitem{Khesin98}
V.~Arnold and B.~Khesin.
\newblock {\em Topological Methods in Hydrodynamics}, volume 125 of {\em
  Applied Mathematical Sciences}.
\newblock Springer-Verlag, New York, 1998.

\bibitem{Benn}
D. Bennequin.
\newblock Entrelacements et \'{e}quations de Pfaff.
\newblock {\em Ast\'{e}risque}, 87--161, Soc. Math. France, Paris, 1983.

\bibitem{Bergery}
L.~Bergery and J. Bourguignon.
\newblock Laplacians and Riemannian submescions with totally geodesic fiberes.
\newblock {\em Illinois J. of Math.}, 26(2):181--199, 1982.

\bibitem{Blair02}
D.~Blair.
\newblock {\em Riemannian geometry of contact and symplectic manifolds}, volume
  203 of {\em Progress in Mathematics}.
\newblock Birkh\"auser Boston Inc., Boston, MA, 2002.

\bibitem{Chandrasekhar57}
S.~Chandrasekhar and P.~Kendall.
\newblock On force-free magnetic fields.
\newblock {\em Astrophys. J.}, 126:457--460, 1957.

\bibitem{Chicone}
C. Chicone.
\newblock The topology of stationary curl parallel solutions of
Euler's equations.
\newblock {\em Israel J. Math.} 39:1-2, 161--166, 1981.


\bibitem{Hamilton85}
S.~S. Chern and R.~S. Hamilton.
\newblock On {R}iemannian metrics adapted to three-dimensional contact
  manifolds.
\newblock In {\em Workshop Bonn 1984 (Bonn, 1984)}, volume 1111 of {\em Lecture
  Notes in Math.}, pages 279--308. Springer, Berlin, 1985.

\bibitem{Eli}
Y. Eliashberg.
\newblock Classification of overtwisted contact structures on $3$-manifolds.
\newblock {\em Invent. Math.} 98:3, 623--637, 1989.

\bibitem{Ghrist00_2}
J.~Etnyre and R.~Ghrist.
\newblock Contact topology and hydrodynamics. {I}. {B}eltrami fields and the
  {S}eifert conjecture.
\newblock {\em Nonlinearity}, 13(2):441--458, 2000.

\bibitem{Giroux01}
E.~Giroux.
\newblock Structures de contact sur les vari\'et\'es fibr\'ees en cercles
  audessus d'une surface.
\newblock {\em Comment. Math. Helv.}, 76(2):218--262, 2001.

\bibitem{Jost02}
J.~Jost.
\newblock Riemannian geometry and geometric analysis.
\newblock {\em Universitext, Berlin: Springer-Verlag, third~ed., 2002}.

\bibitem{Komendarczyk_nodal}
R.~Komendarczyk.
\newblock On the contact geometry of nodal sets, {\tt math.dg/0402070}.
\newblock To appear, {\em  Trans. of Amer. Math. Soc.}, 2004.

\bibitem{Nicolaescu98}
L.~Nicolaescu.
\newblock Adiabatic limits of the {S}eiberg-{W}itten equations on {S}eifert
  manifolds.
\newblock {\em Comm. Anal. Geom.}, 6(2):331--392, 1998.

\bibitem{Uhlenbeck76}
K.~Uhlenbeck.
\newblock Generic properties of eigenfunctions.
\newblock {\em Amer. J. Math.}, 98(4):1059--1078, 1976.

\bibitem{Yoshida90}
Z.~Yoshida and Y.~Giga.
\newblock Remarks on spectra of operator rot.
\newblock {\em Math. Z.}, 204(2):235--245, 1990.

\end{thebibliography}

\bibliographystyle{plain}

\end{document}